\documentclass[10pt]{article}
\usepackage{pb-diagram}
\usepackage{amsmath}
\usepackage{amssymb}
\usepackage{amscd}

\newtheorem{Df}{Definition}[section]
\newtheorem{Te}[Df]{Theorem}
\newtheorem{Po}[Df]{Proposition}
\newtheorem{Cr}[Df]{Corollary}
\newtheorem{Lm}[Df]{Lemma}
\newtheorem{Ca}[Df]{Claim}
\newtheorem{Cn}[Df]{Conjecture}
\newtheorem{Ex}[Df]{Example}

\newcommand{\Bdf}{\begin{Df}}
\newcommand{\Edf}{\end{Df}}
\newcommand{\Bte}{\begin{Te}}
\newcommand{\Ete}{\end{Te}}
\newcommand{\Bpo}{\begin{Po}}
\newcommand{\Epo}{\end{Po}}
\newcommand{\Bcr}{\begin{Cr}}
\newcommand{\Ecr}{\end{Cr}}
\newcommand{\Blm}{\begin{Lm}}
\newcommand{\Elm}{\end{Lm}}
\newcommand{\Bca}{\begin{Ca}}
\newcommand{\Eca}{\end{Ca}}
\newcommand{\Bcn}{\begin{Cn}}
\newcommand{\Ecn}{\end{Cn}}
\newcommand{\Bex}{\begin{Ex}}
\newcommand{\Eex}{\end{Ex}}
\newcommand{\Bdm}{{\it Proof.}\ }
\newcommand{\Edm}{\rule{2mm}{2mm}}

\begin{document}

\title{\bf{A criterion for homogeneous potentials to be 3-Calabi-Yau}}
\author{\bf{Roland Berger and Andrea Solotar\footnote{This work has been partially supported by the projects UBACYT 20020100100475, PIP-CONICET 
112-200801-00487, PICT-2007-02182 and MATHAMSUD-NOCOSETA. A.Solotar is a research member of CONICET (Argentina).}}}
\date{}
\maketitle
\begin{abstract}
Among the homogeneous potentials $w$ of degree $N+1$ in $n$ variables, it is an open problem to find precisely which of the $w$'s are 3-Calabi-Yau, 
although several examples are known. In this paper, we give a necessary and sufficient condition for this to hold when the algebra $A$ 
defined by the potential $w$ is $N$-Koszul of global dimension 3. As an application, we study skew polynomial algebras over non-commutative quadrics and 
we recover two families of 3-Calabi-Yau potentials which have recently appeared in the literature.
\end{abstract} 

\section{Introduction}
Following Ginzburg~\cite{vg:cy}, $d$-Calabi-Yau algebras $A$ are defined by some natural finiteness constraints and a duality condition involving 
the Hochschild cohomology $\mbox{H}^{\bullet}(A, A\otimes A^{op})$. After examination of various examples, several authors conjectured that any 3-Calabi-Yau 
algebra $A$ (satisfying some more or less natural assumptions) can be defined nicely by generators and relations from a certain non-commutative 
polynomial depending on $A$ and called \emph{potential}. In this situation, one says that $A$ is a potential algebra or that $A$ is derived from a potential 
(see the precise definition in Section 2 below). The best result in this direction has been obtained by Van den Bergh, who proved that any complete 
3-Calabi-Yau algebra is derived from a potential~\cite{vdb:superpo}. A proof for graded algebras $A$ generated in degree one had been previously given 
by Bocklandt~\cite{bock:graded}.

Actually, a potential algebra is derived in a standard way from any potential, i.e. from any non-commutative polynomial.
However it is not known which are the potentials $w$ such that the algebra $A$ derived from $w$ is 3-Calabi-Yau, even in the graded case. In this paper, 
we restrict ourselves to \emph{homogeneous} potentials $w$ in $n\geq 1$ variables of degree one. Denote by $N+1$ the degree of $w$, where $N\geq 2$, so that 
the graded algebra $A$ derived from $w$ is $N$-homogeneous. It is known  
(Proposition 5.2 in~\cite{rbnm:kogo}) that if $A$ is AS-Gorenstein of global dimension 3 (in particular if $A$ is 3-Calabi-Yau~\cite{rbrt:pbw}), 
then $A$ is $N$-Koszul with $n\geq 2$, and its Hilbert series is the following:
\begin{equation} \label{hs}
h_A(t)=(1-nt + nt^N - t^{N+1})^{-1}.
\end{equation}
The main result of this paper is the following theorem proved in Section 2. In all the paper, 
the basic field $k$ has characteristic zero.
\Bte \label{main}
For any $n\geq 2$ and $N\geq 2$, let $A$ be a graded algebra in $n$ generators $x_1, \ldots , x_n$, derived from a homogeneous potential 
$w$ of degree $N+1$. Assume that $A$ is $N$-Koszul. Then $A$ is 3-Calabi-Yau if and only if the Hilbert series of $A$ is given by (\ref{hs}). 
\Ete

Denoting by $V$ and $R$ the spaces of generators and relations of such an $N$-Koszul potential algebra $A$, condition (\ref{hs}) is equivalent to saying 
that the three following facts hold:

(i) $A$ has global dimension 3,

(ii) the cyclic partial derivatives $\partial_{x_1}(w), \ldots , \partial_{x_n}(w)$ are $k$-linearly independent,

(iii) the cyclic sum $c(w)$ of $w$ generates the space $R_{N+1}=(R\otimes V)\cap (V\otimes R)$.
\\ 

In Theorem 6.8 of~\cite{bsw:higher}, Bocklandt, Schedler and Wemyss have proved that if $A$ is an algebra defined by a twisted potential $w$ and $A$ 
is $N$-Koszul, then $A$ is twisted $d$-Calabi-Yau if and only if a certain complex defined from $w$ $-$ a bimodule version of a complex previously 
considered by Dubois-Violette~\cite{dv:multi} $-$ is exact. In the non-twisted case with $d=3$, our Theorem \ref{main} is an 
improvement of their result: in order to conclude that $A$ is 3-Calabi-Yau, it is sufficient to know the Hilbert series, while 
exactness is in general hard to prove. It would be interesting to extend Theorem \ref{main} to the general setting of~\cite{bsw:higher}.

The examples of applications of Theorem \ref{main} we have in mind are quadratic (i.e., for $N=2$) and they have recently appeared 
in various contexts. The example due to Smith is constructed from the octonions~\cite{smith:octo}, and more generally those due to 
Su\'arez-Alvarez are constructed from the oriented Steiner triple systems~\cite{msa:steiner}. The examples due to Berger and Pichereau come from 
a way to embed any non-degenerate non-commutative quadric (not necessarily 3-Calabi-Yau) into a 3-Calabi-Yau potential algebra by adding a 
variable~\cite{rbap:cyp}. Actually, as shown in the cited articles, each of these examples is a skew polynomial algebra $A$ 
over a non-degenerate non-commutative quadric $\Gamma$, in other words $A$ is an Ore extension of $\Gamma$. Then, using the basic properties of $\Gamma$ 
as stated in~\cite{rb:gera} (see also~\cite{rbap:cyp}), 
the 3-Calabi-Yau property is immediate from the following consequence (proved in Section 4 below) of Theorem \ref{main}.
\Bcr \label{cormain}
For any $n\geq 2$, let $\Gamma$ be a non-degenerate non-commutative quadric in $n$ variables $x_1, \ldots, x_n$ of degree 1. Let $z$ be an extra variable of 
degree 1. Let $A$ be an algebra defined by a non-zero cubic potential $w$ in the variables $x_1, \ldots , x_n,\, z$. Assume that the graded algebra $A$ is 
isomorphic to a skew polynomial algebra $\Gamma[z;\, \sigma;\, \delta]$ over $\Gamma$ in the variable $z$, defined by a 0-degree homogeneous automorphism 
$\sigma$ of $\Gamma$ and a 1-degree homogeneous $\sigma$-derivation $\delta$ of $\Gamma$. Then $A$ is Koszul and 3-Calabi-Yau.
\Ecr

In the situation considered in~\cite{rbap:cyp}, where the potential $w$ is defined by $w=uz$ for $u$ the relation of the quadric $\Gamma$, Ulrich Kr\"{a}hmer has asked the first 
author whether the automorphism $\sigma$ is related to the automorphism of Van den Bergh's duality of $\Gamma$. We answer this question by 
proving in Section 4 that the dualizing bimodule of $\Gamma$ is isomorphic to $_{\sigma ^{-1}}\Gamma$. 
\\ 

\emph{Acknowledgements.} We thank Thierry Lambre for fruitful discussions during the preparation of this work. 
We are grateful to Ulrich Kr\"{a}hmer for his very interesting question and to Eduardo Marcos for a careful reading of this article. The second author thanks the University of Saint-Etienne for her stays in the 
Faculty of Sciences.

\setcounter{equation}{0}

\section{3-Calabi-Yau potential algebras}

Throughout the article, $k$ will denote a field of characteristic zero and $V$ a $k$-vector space of finite dimension $n\geq 1$. We fix a basis 
$X=\{x_1, \ldots, x_n\}$ of $V$. The tensor algebra of $V$ is denoted by $T(V)$ and $T(V)^e=T(V)\otimes T(V)^{op}$. The symbol $\otimes_k$ will always be 
denoted by $\otimes$. The basis of $T(V)$ consisting of the non-commutative monomials in 
$x_1, \ldots , x_n$ is denoted by $\langle X \rangle$. The subspace of $T(V)$ generated by the commutators is denoted by $[T(V), T(V))]$. The elements 
of the vector space $$Pot(V)=T(V)/[T(V), T(V)]$$ also denoted by $Pot(x_1, \ldots , x_n)$, are called \emph{potentials} of $V$ or potentials in 
the variables $x_1, \ldots ,x_n$. The algebra $T(V)$ is graded by $\deg(x_i)=1$ for all $i$. Since $[T(V),T(V)]$ is homogeneous, the space $Pot(V)$ 
inherits the grading. The linear map $c:T(V) \rightarrow T(V)$ defined on monomials $a=a_1\ldots a_r$ of degree $r$ when all the $a_i$'s are in $X$, 
as the cyclic sum 
$$c(a) = \sum_{1\leq i \leq r} a_i\ldots a_r a_1 \ldots a_{i-1}$$
induces $\tilde{c}:Pot(V) \rightarrow T(V)$ which actually does not depend on the choice of the basis $X$. Since the characteristic of $k$ is zero, 
this map $\tilde{c}$ defines a linear isomorphism from $Pot(V)$ onto $Im(c)$.

For any $x\in X$, the cyclic derivative $\partial_x : Pot(V) \rightarrow T(V)$ is the linear map defined on any $p\in\langle X \rangle $ by 
$$\partial_x(p)= \sum_{p=uxv} vu,$$
where $u$ and $v$ are in $\langle X \rangle $. The ``ordinary'' partial derivative $\frac{\partial}{\partial x} : T(V)
\rightarrow T(V)\otimes T(V)$ (see~\cite{vdb:ober}) is the linear map defined on any monomial $p$ by 
$$\frac{\partial p}{\partial x} = \sum_{p=uxv} u\otimes v,$$
which will be written as 
\begin{equation} \label{symb}
\frac{\partial p}{\partial x} = \sum_{1,2} \left(\frac{\partial p}{\partial x}\right)_1
\otimes \left(\frac{\partial p}{\partial x}\right)_2 \ .
\end{equation}
Considering $T(V)$ as the natural $T(V)$-bimodule, hence as a right $T(V)^e$-module, we verify that 
$$1 \cdot \frac{\partial p}{\partial x}=\partial_x(p).$$
Finally, for $x$ and $y$ in $X$, the second partial derivative $\frac{\partial^2 }{\partial x \partial y} : Pot(V)\rightarrow T(V)\otimes T(V)$ 
is defined by $$\frac{\partial^2 }{\partial x \partial y}=\frac{\partial}{\partial x} \circ \partial_{y}.$$ 

For any $w\in Pot(V)$, we recall two basic formulas concerning the previous derivatives which can be easily proved. The first one is the non-commutative
Euler relation (generalizing the well-known Euler relation for homogeneous commutative polynomials):
\begin{equation} \label{eul}
\sum_{1\leq i \leq n} \partial_{x_i} (w) x_i = \sum_{1\leq i \leq n} x_i \partial_{x_i} (w) = c(w),
\end{equation}
where the constant term of $w$ is assumed to be zero.
The second one is the symmetry of the non-commutative Hessian~\cite{vdb:ober}:
\begin{equation} \label{hess}
\tau \left(\frac{\partial^2 w}{\partial x \partial y}\right)= \frac{\partial^2 w}{\partial
y \partial x}
\end{equation} 
where $\tau :T(V)\otimes T(V) \rightarrow T(V)\otimes T(V)$ is the flip $a\otimes b \mapsto b\otimes a$.
\Bdf
For any $w\in Pot(V)$, let $I(\partial_x(w);\,x\in X)$ denote the two-sided ideal generated by all the cyclic partial derivatives of $w$. We say that 
the associative $k$-algebra 
$$A=A(w)=T(V)/I(\partial_x(w);\,x\in X)$$
is derived from the potential $w$, or that it is the potential algebra defined from $w$.
\Edf

For the rest of this section, let us fix an integer $N\geq 2$ and a non-zero homogeneous potential $w$ of degree $N+1$. The space of homogeneous 
potentials of degree $N+1$ is the following
$$Pot(V)_{N+1}=\frac{V^{\otimes (N+1)}}{\sum_{i+j=N+1} [V^{\otimes i} ,\, V^{\otimes j}]} = 
\frac{V^{\otimes (N+1)}}{\sum_{i+j=N-1} V^{\otimes i} \otimes [V,V] \otimes V^{\otimes j}}\,,$$
and the class $w\in Pot(V)_{N+1}$ will be often defined by a representative denoted again by $w$. The $k$-algebra $A$ derived 
from $w$ is $\mathbb{N}$-graded and the relations $\partial_x(w)$, $x\in X$, are homogeneous of degree $N$. So the graded algebra $A$ 
is $N$-homogeneous~\cite{bdvw:homog}. Let us denote by $R$ the space of relations of $A$, i.e. the subspace of $V^{\otimes N}$ generated by the relations 
$\partial_x(w)$, $x\in X$. The subspace $R_{N+1}=(R\otimes V)\cap (V\otimes R)$ of $V^{\otimes (N+1)}$ appears 
in the Koszul complex of $A$~\cite{rb:nonquad} and moreover, the non-commutative Euler relation (\ref{eul}) shows that 
\begin{equation} \label{fond}
c(w) \in R_{N+1}.
\end{equation}
Remark that $c(w) \neq 0$ since $w\neq 0$. The element $c(w)$ will play the role of the volume form in the non-commutative setting. 
For convenience, we will sometimes omit the unadorned symbols $\otimes$, e.g. we will write $R_{N+1}=RV\cap VR$, $A\otimes A = AA, \ldots$

Recall that the bimodule Koszul complex of $A$ starts as follows~\cite{rbnm:kogo}:
\begin{equation} \label{start1}
A R A \stackrel{d_2}\longrightarrow
A V\! A \stackrel{d_1}\longrightarrow A A \rightarrow 0,
\end{equation}
where $d_1$ and $d_2$ are $A$-$A$-linear and their restrictions to $V$ (resp. to $R$) are defined by $$d_1(v)= v\otimes 1 -1 \otimes v \in AA,$$ 
$$d_2 (v_1 \ldots v_N)= \sum_{1\leq i \leq N} (v_1\ldots v_{i-1})\otimes v_i\otimes (v_{i+1}\ldots v_N)\in AV\!A,$$ 
for any $v, v_1, \ldots v_N$ in $V$. Moreover, via the multiplication $\mu : AA \rightarrow A$, the sequence (\ref{start1}) can be extended to a minimal 
projective resolution of $A$ in the category $A$-grMod-$A$ of graded $A$-bimodules. 

The graded vector space $\ker d_2$ lives in degrees $\geq N+1$ and the linear map $\varphi: R_{N+1} \rightarrow ARA$ defined by 
$$\varphi\,(\sum_{1\leq i \leq n} x_i \otimes u_i = \sum_{1\leq i \leq n} v_i \otimes x_i\,)= \sum_{1\leq i \leq n}x_i\otimes u_i \otimes 1 - 1 \otimes v_i 
\otimes x_i$$
where $u_i$ and $v_i$ are in $R$, is an isomorphism from $R_{N+1}$ to $(\ker d_2)_{N+1}$. Choose a graded subspace $E$ of $\ker d_2$ such that 
$$\ker d_2= E\oplus \sum_{i+j \geq 1} A_i (\ker d_2)_{N+1} A_j.$$
Note that $E$ lives in degrees $\geq N+1$ and that $E_{N+1}= (\ker d_2)_{N+1}$. Then the complex 
\begin{equation} \label{start2}
AEA \stackrel{d_3}\longrightarrow A R A \stackrel{d_2}\longrightarrow
A V\! A \stackrel{d_1}\longrightarrow A A \rightarrow 0,
\end{equation}
where $d_3$ is the $A$-$A$-linear extension of the inclusion of $E$ into $ARA$, can be extended via the multiplication $\mu$ of $A$ 
to a minimal projective resolution of $A$ in $A$-grMod-$A$. Actually, using the isomorphism $\varphi$, we will assume that $E_{N+1}=R_{N+1}$ and 
that $d_3$ coincides with $\varphi$ on $R_{N+1}$. So $E_{N+1}$ contains the non-zero element $c(w)$. 
\Blm \label{smalm}
Keep the above notation and assumptions. The global dimension of $A$ is equal to 3 if and only if $d_3$ is injective. In this case, $A$ is $N$-Koszul if and only if $E=R_{N+1}$. 
\Elm
\Bdm
The global dimension of $A$ is the length of a minimal projective resolution, hence the first equivalence follows. The second one is clear from 
the definition of $N$-Koszul algebra~\cite{rb:nonquad}. 
\Edm 
\\

The next result is an immediate consequence of Proposition 5.2 in~\cite{rbnm:kogo}. For the convenience of the reader, we give here a self-contained proof.
\Bpo \label{Gore}
Let $k$ be a field of characteristic zero. Let $V\neq 0$ be an $n$-dimensional vector space. Let $w$ be a non-zero homogeneous potential of $V$ of 
degree $N+1$ with $N\geq 2$. Let $A=A(w)$ be the potential algebra defined by $w$. 
Denote by $R$ the space of relations of $A$. Assume that $A$ is AS-Gorenstein of global dimension 3. Then 

(i) $\dim R =n$ (with $n\geq 2$), so that $(\partial_{x_i} (w))_{1\leq i \leq n}$ is a basis of $R$,

(ii) $E=R_{N+1}$ which is one-dimensional generated by $c(w)$,

(iii) $A$ is $N$-Koszul,

(iv) the Hilbert series of the graded algebra $A$ is given by $h_A(t)=(1-nt + nt^N - t^{N+1})^{-1}$.
\Epo
\Bdm
Since the global dimension of $A$ is equal to 3, the map $d_3$ in (\ref{start2}) is injective, so 
$$0\rightarrow AEA \stackrel{d_3}\longrightarrow A R A \stackrel{d_2}\longrightarrow
A V\! A \stackrel{d_1}\longrightarrow A A \rightarrow 0$$
is a minimal projective resolution of $A$ via $\mu$. Thus the AS-Gorenstein symmetry shows that $\dim R = \dim V$ and $\dim E=1$, hence we get (i) 
(where $n=1$ is easily ruled out) and (ii). 
Next we use the Lemma. The Hilbert series is immediately obtained from the exact complex 
$$0\rightarrow AEA \stackrel{d_3}\longrightarrow A R A \stackrel{d_2}\longrightarrow
A V\! A \stackrel{d_1}\longrightarrow A A \stackrel{\mu}\longrightarrow A \rightarrow 0. \ \Edm$$

Throughout the rest of this section, $A$ is derived from a non-zero homogeneous potential $w$ of degree $N+1$ ($N\geq 2$) in $n\geq 1$ variables, 
and we use the previous notation. Before proving Theorem \ref{main} of the Introduction, we obtain some general results without assuming yet that $h_A(t)=(1-nt + nt^N - t^{N+1})^{-1}$. In particular, we just have $\dim R \leq n$. We want to study the self-duality of the following complex denoted by $C_w$:
\begin{equation} \label{sdc}
0 \rightarrow Akc(w)A \stackrel{d_3}\longrightarrow A R A \stackrel{d_2}\longrightarrow
A V\! A \stackrel{d_1}\longrightarrow AkA \rightarrow 0,
\end{equation}
where $k$ in $AA=AkA$ is the subspace of $AA$ generated by $1\otimes 1$, and $kc(w)$ in $Akc(w)A$ is the subspace of $AR_{N+1}A$ generated 
by $c(w)$. Actually, (\ref{sdc}) is a subcomplex of 
$$0\rightarrow AEA \stackrel{d_3}\longrightarrow A R A \stackrel{d_2}\longrightarrow
A V\! A \stackrel{d_1}\longrightarrow AkA \rightarrow 0.$$
We already know that it is exact at $AVA$ and that its homology at $AkA$ 
is isomorphic to $A$ (using the multiplication $\mu : AkA \rightarrow A$). Let $r_i= \partial_{x_i}(w)$ for $1\leq i \leq n$. According to (\ref{eul}), 
we have 
\begin{equation} \label{eul2}
d_3(c(w))= \sum_{1\leq i \leq n} x_i \otimes r_i \otimes 1 - 1 \otimes r_i \otimes x_i \,.
\end{equation}

Next we recall some general facts about duality of bimodules: let $A_1$ and $A_2$ be associative $k$-algebras and let $M$ be an $A_1$-$A_2$-bimodule. Denote by 
$A_1 \stackrel{o}\otimes A_2$ (resp. $A_1 \stackrel{i}\otimes A_2$) the $A_1$-$A_2$ (resp. $A_2$-$A_1$) bimodule $A_1\otimes A_2$ for the outer 
(resp. inner) action. Set 
$$M^{\vee} = \mbox{Hom}_{A_1-A_2} (M, A_1 \stackrel{o}\otimes A_2),$$
so that $M^{\vee}$ is an $A_2$-$A_1$-bimodule whose action comes from $A_1 \stackrel{i}\otimes A_2$. For any finite-dimensional $k$-vector space  $E$ whose dual 
will be denoted by $E^{\ast}$, we have a natural isomorphism 
$$\theta : (A_1EA_2)^{\vee} \rightarrow A_2 E^{\ast} A_1$$
of $A_2$-$A_1$-bimodules that we are going to describe. Firstly, we have naturally 
$$\mbox{Hom} _{A_1-A_2} (A_1EA_2, A_1 A_2) \cong \mbox{Hom}_k (E, A_1 A_2).$$
Secondly, if $\gamma \in \mbox{Hom} _k (E, A_1 A_2)$, we define $\tilde{\gamma}\in E^{\ast}A_1A_2$ by
$$\tilde{\gamma}= \sum_{i\in I} v_i^{\ast} \otimes \gamma (v_i),$$
where $(v_i)_{i\in I}$ is a basis of $E$ and $(v_i^{\ast})_{i\in I}$ is its dual basis. Then 
$$\mbox{Hom}_k (E, A_1A_2) \cong E^{\ast} A_1A_2$$
by the $k$-linear isomorphism $\gamma \mapsto \tilde{\gamma}$ whose inverse isomorphism is given by 
$$\phi a_1 a_2 \mapsto (v\mapsto \phi (v) a_1a_2),\ \  \phi \in E^{\ast}, \ a_1 \in A_1,\ a_2 \in A_2.$$
Finally $E^{\ast}A_1A_2 \cong A_2 E^{\ast} A_1$ is natural with respect to the obvious $A_2$-$A_1$-bimodule structures. Composing all these maps, 
we get the isomorphism $\theta$. 

Throughout the sequel, we take $A_1=A_2=A$.
For any chain complex $(C,d)$ of $A$-bimodules, the dual complex 
$$C^{\vee}= \mbox{Hom}_{A-A} (C, A \stackrel{o}\otimes A)$$
(as usual, Hom is graded) is a chain complex of $A$-bimodules whose differential 
$$d_{1-n}^{\vee} : C_n^{\vee}= \mbox{Hom}_{A-A} (C_{-n}, AA) \rightarrow C_{n-1}^{\vee}$$ 
is defined for any $n\in \mathbb{Z}$ by 
\begin{equation} \label{duald}
d_{1-n}^{\vee}(f)= - (-1)^n f\circ d_{1-n},
\end{equation}
where $d_{1-n} : C_{1-n} \rightarrow C_{-n}$ and $f : C_{-n} \rightarrow AA$. Note that the sign $(-1)^n$ in this definition comes from the Koszul rule 
(see for example~\cite{nb:hom}, formula (1), p. 81). 

Let us compute the differential of the dual complex $C_w^{\vee}$:
\begin{equation} \label{dsdc}
0 \rightarrow Ak^{\ast}A \stackrel{d_1^{\ast}}\longrightarrow A V^{\ast} A \stackrel{d_2^{\ast}}\longrightarrow
A R^{\ast} A \stackrel{d_3^{\ast}}\longrightarrow Akc(w)^{\ast}A \rightarrow 0,
\end{equation}
where $d_i^{\ast}$ denotes the image of $d_i^{\vee}$ via the isomorphism $\theta$. Firstly, we have $d_1^{\ast}=\theta \circ d_1^{\vee}\circ \theta^{-1}$. 
From the dual basis $1^{\ast}\in k^{\ast}$ of $1\in k$ , we get that $\theta^{-1}(1^{\ast})$ coincides with the identity map $1_{AA}$ of $AA$, thus 
$\gamma=d_1^{\vee}(\theta^{-1}(1^{\ast}))$ is defined in $V$ by $v\mapsto -1_{AA}\circ d_1(v)=1\otimes v -v\otimes 1$ (note that $n=0$ in (\ref{duald})). 
We get $\tilde{\gamma} =
\sum_{1\leq i \leq n} x_i^{\ast} \otimes (1\otimes x_i -x_i \otimes 1)$, where $(x_i^{\ast})_{1\leq i \leq n}$ is the dual basis of the basis 
$(x_i)_{1\leq i \leq n}$ of $V$. We conclude that 
\begin{equation} \label{d1ast}
d_1^{\ast}(1^{\ast})= \sum_{1\leq i \leq n} x_i \otimes x_i^{\ast} \otimes 1-1 \otimes x_i^{\ast} \otimes x_i.
\end{equation}

Secondly, we have $d_2^{\ast}=\theta \circ d_2^{\vee}\circ \theta^{-1}$. Recall that $r_1= \partial_{x_1} w, \ldots , r_n=\partial_{x_n} w$ generate $R$. 
Fix a part $J$ of $\{1, \ldots , n\}$ such that $(r_j)_{j\in J}$ is a basis of $R$, and denote by $(r_j^{\ast})_{j\in J}$ its dual basis. We get that 
$\theta^{-1}(x_i^{\ast})$ coincides with the map $AV\!A \rightarrow AA$, $ax_j b \mapsto \delta_{ij} ab$, thus $\gamma_i=d_2^{\vee}(\theta^{-1}(x_i^{\ast}))$ 
is defined on $R$ by $r_j \mapsto \theta^{-1}(x_i^{\ast}) \circ d_2 (r_j)$ since $n=-1$ in (\ref{duald}). Using the symbolic writing (\ref{symb}), we have 
\begin{equation} \label{d2}
d_2(r_j)=\sum_{1\leq s\leq n} \sum_{1,2} \left(\frac{\partial r_j}{\partial x_s}\right)_1
\otimes x_s \otimes \left(\frac{\partial r_j}{\partial x_s}\right)_2 \ .
\end{equation}
Therefore $\gamma_i(r_j)= \sum_{1,2} \left(\frac{\partial r_j}{\partial x_i}\right)_1 \otimes \left(\frac{\partial r_j}{\partial x_i}\right)_2$, which 
together with the equality $\tilde{\gamma_i}= \sum_{j\in J} r_j^{\ast}\otimes \gamma_i(r_j)$ 
implies that, for any $1\leq i \leq n$, 
\begin{equation} \label{d2ast}
d_2^{\ast}(x_i^{\ast})= \sum_{j \in J} \sum_{1,2} \left(\frac{\partial r_j}{\partial x_i}\right)_2
\otimes r_j^{\ast} \otimes \left(\frac{\partial r_j}{\partial x_i}\right)_1.
\end{equation}

Thirdly, we have $d_3^{\ast}=\theta \circ d_3^{\vee}\circ \theta^{-1}$. We get that $\theta^{-1}(r_i^{\ast})$ (for $i\in J$) coincides with 
the map $ARA \rightarrow AA$, $ar_j b \mapsto \delta_{ij} ab$ with $j\in J$, thus $\gamma_i=d_3^{\vee}(\theta^{-1}(r_i^{\ast}))$ is defined on $kc(w)$ by 
$c(w) \mapsto -\theta^{-1}(r_i^{\ast}) \circ d_3 (c(w))$ since $n=-2$ in (\ref{duald}). Using (\ref{eul2}), we obtain for any $i \in J$, 
$$\gamma_ i(c(w))=1\otimes x_i - x_i \otimes 1 +\sum_{j\notin J} r_i^{\ast}(r_j)(1\otimes x_j - x_j \otimes 1).$$
Thus, for any $i \in J$, we conclude that 
\begin{equation} \label{d3ast}
d_3^{\ast}(r_i^{\ast})=  x_i \otimes c(w)^{\ast} \otimes 1 - 1 \otimes c(w)^{\ast} \otimes x_i + 
\sum_{j\notin J} r_i^{\ast}(r_j)(x_j \otimes c(w)^{\ast} \otimes 1 - 1 \otimes c(w)^{\ast} \otimes x_j).
\end{equation}
In particular, if the $A$-$A$-linear map $\mu_w: A(kc(w)^{\ast}) A \rightarrow A$ is defined by $\mu_w(c(w)^{\ast})=1$, then we form the augmented complex 
$$C_w^{\vee} \stackrel{\mu_w}{\longrightarrow} A \rightarrow 0.$$
Moreover, the $A$-$A$-linear map $f_0 : AkA \rightarrow A(kc(w)^{\ast}) A$ defined by $f_0(1)=c(w)^{\ast}$ is such that $\mu_w \circ f_0 = \mu$.

Consider the diagram 
\begin{eqnarray} \label{self}
0 \rightarrow A(kc(w))A \stackrel{d_{3}}{\longrightarrow} & ARA
\stackrel{d_{2}}{\longrightarrow} 
& AV\! A \stackrel{d_{1}}{\longrightarrow} AkA \longrightarrow 0
\nonumber \\
f_{3} \downarrow \ \ \ \ \ \ \  & f_{2} \downarrow \ \ \ \ \ \ \ &  f_{1}
\downarrow \ \ \ \ \ \ \ \ \  f_{0} \downarrow  \\
0 \longrightarrow Ak^{\ast}A \stackrel{d_{1}^{\ast}}{\longrightarrow} & AV^{\ast}A \stackrel{d_{2}^{\ast}}{\longrightarrow}
& AR^{\ast}A\stackrel{d_{3}^{\ast}}{\longrightarrow}  A(kc(w)^{\ast}) A \rightarrow 0
\nonumber 
\end{eqnarray}
where the $A$-$A$-linear maps $f_1$, $f_2$ and $f_3$ are given by
\begin{equation}
f_1(x_i)= \left\{ \begin{array}{ll}
r_i^{\ast} & \mbox{if}\ i\in J\\ 
0 & \mbox{if}\ i\notin J
\end{array}
\right.
,\ \ f_2(r_i)=x_i^{\ast}\ \mbox{for\ any}\ i\in J, \ \ f_3(c(w))= 1^{\ast}.
\end{equation}
Recall that $J$ is a part of $\{1,\ldots, n\}$ such that $(r_j)_{j\in J}$ is a basis of $R$ and $(r_j^{\ast})_{j\in J}$ is its dual basis. 
Clearly, $f_0$ and $f_3$ are bijective, $f_1$ is surjective 
and $f_2$ is injective. The following proposition shows that the diagram (\ref{self}) is commutative if and only if $\dim R =n$. 
\Bpo  \label{commut}
(i) The central square in (\ref{self}), i.e. the square limited by $f_1$ and $f_2$, is commutative.

(ii) Each of both remaining squares is commutative if and only if $\dim R = n$.

(iii) If $C_w^{\vee}$ is exact at $AV^{\ast}A$, i.e. if $\emph{H}^1 (A, A\stackrel{o}\otimes A)=0$, then $\dim R = n$.
\Epo
\Bdm
(i) For any $i \in J$, we have
$$f_1\circ d_2 (r_i)= f_1 (\sum_{1\leq j\leq n} \sum_{1,2} \left(\frac{\partial r_i}{\partial x_j}\right)_1
\otimes x_j \otimes \left(\frac{\partial r_i}{\partial x_j}\right)_2)= \sum_{j \in J} \sum_{1,2} \left(\frac{\partial r_i}{\partial x_j}\right)_1
\otimes r_j^{\ast} \otimes \left(\frac{\partial r_i}{\partial x_j}\right)_2$$
$$d_2^{\ast} \circ f_2 (r_i)=d_2^{\ast}(x_i^{\ast})= \sum_{j \in J} \sum_{1,2} \left(\frac{\partial r_j}{\partial x_i}\right)_2
\otimes r_j^{\ast} \otimes \left(\frac{\partial r_j}{\partial x_i}\right)_1.$$
So, in these both sums, the coefficients of $r_j^{\ast}$ are respectively 
$$\frac{\partial }{\partial x_j}\circ \partial_{x_i} (w)= \frac{\partial^2 w}{\partial x_j \partial x_i},$$
$$\tau \left(\frac{\partial }{\partial x_i}\circ \partial_{x_j} (w)\right)= \tau \left(\frac{\partial^2 w}{\partial x_i \partial x_j}\right).$$
Thus they are equal by the symmetry (\ref{hess}) of the non-commutative Hessian.

(ii) On one hand, from $c(w)=\sum_{1\leq i \leq n} x_i\otimes r_i \otimes 1 - 1 \otimes r_i \otimes x_i$, we get 
$$f_2\circ d_3 (c(w))= \sum_{i \in J} x_i \otimes x_i^{\ast} \otimes 1 - 1 \otimes x_i^{\ast} \otimes x_i + 
\sum_{i\notin J, j\in J} r_j^{\ast}(r_i)(x_i \otimes x_j^{\ast} \otimes 1 - 1 \otimes x_j^{\ast} \otimes x_i).$$
From (\ref{d1ast}), we have
$$d_1^{\ast}\circ f_3 (c(w))= d_1^{\ast}(1^{\ast}) = \sum_{i \in J} x_i \otimes x_i^{\ast} \otimes 1 - 1 \otimes x_i^{\ast} \otimes x_i 
+ \sum_{i \notin J} x_i \otimes x_i^{\ast} \otimes 1 - 1 \otimes x_i^{\ast} \otimes x_i.$$
If $\dim R = n$, then $f_2\circ d_3=d_1^{\ast}\circ f_3$. The converse comes from the linear independence of $x_1^{\ast}, \ldots, x_n^{\ast}$.

On the other hand, we have for any $1\leq i \leq n$,
$$f_0 \circ d_1 (x_i) = x_i \otimes c(w)^{\ast} \otimes 1 - 1 \otimes c(w)^{\ast} \otimes x_i.$$
If $i\notin J$, then $d_3^{\ast}\circ f_1 (x_i)= 0$. If $i\in J$, we have 
$$d_3^{\ast}\circ f_1 (x_i)=d_3^{\ast}(r_i^{\ast})=x_i \otimes c(w)^{\ast} \otimes 1 - 1 \otimes c(w)^{\ast} \otimes x_i + \sum_{j\notin J} r_i^{\ast}(r_j)(x_j \otimes c(w)^{\ast} \otimes 1 - 1 \otimes c(w)^{\ast} \otimes x_j).$$
If $\dim R = n$, then $f_0\circ d_1=d_3^{\ast}\circ f_1$. The converse comes from $c(w)^{\ast}\neq 0$.

(iii) Assume that $C_w^{\vee}$ is exact at $AV^{\ast}A$ and that there exists $i\notin J$. Set $r_i=\sum_{j \in J}\lambda_{ij} r_j$ where 
$\lambda_{ij} \in k$, and apply $\frac{\partial}{\partial x_s}$ to this equality for $1\leq s\leq n$. We obtain 
\begin{equation} \label{inter}
\sum_{1,2} \left(\frac{\partial r_i}{\partial x_s}\right)_1
\otimes  \left(\frac{\partial r_i}{\partial x_s}\right)_2 = \sum_{j \in J} \sum_{1,2} \lambda_{ij} \left(\frac{\partial r_j}{\partial x_s}\right)_1
\otimes  \left(\frac{\partial r_j}{\partial x_s}\right)_2 .
\end{equation}
However (\ref{d2ast}) shows that $d_2^{\ast} (x_i^{\ast}-\sum_{j\in J} \lambda_{ij} x_j^{\ast})$ is equal to
$$\sum_{1\leq s\leq n} \left(\sum_{1,2} \left(\frac{\partial r_s}{\partial x_i}\right)_2
\otimes r_s^{\ast} \otimes \left(\frac{\partial r_s}{\partial x_i}\right)_1 - \sum_{j \in J} \sum_{1,2} \lambda_{ij} \left(\frac{\partial r_s}
{\partial x_j}\right)_2 \otimes r_s^{\ast} \otimes \left(\frac{\partial r_s}{\partial x_j}\right)_1\right).$$
For each $s$, the terms with $r_s^{\ast}$ vanish by applying the Hessian symmetry (\ref{hess}) to (\ref{inter}). Hence 
$x_i^{\ast}-\sum_{j\in J} \lambda_{ij} x_j^{\ast}$ is a 1-cocycle of $C_w^{\vee}$. Since the $d_i$'s are homogeneous of degree 0, the same holds 
for the $d_i^{\ast}$'s, where $k^{\ast}$, $V^{\ast}$, $R^{\ast}$, $kc(w)^{\ast}$ are respectively concentrated in degrees 0, $-1$, $-N$ and $-N-1$. 
In particular, the 1-coboundaries live in degrees $\geq 0$. Thus the exactness of $C_w^{\vee}$ at $AV^{\ast}A$ implies that 
$x_i^{\ast}-\sum_{j\in J} \lambda_{ij} x_j^{\ast}=0$, which is a contradiction. 
\Edm
\\ 

It is possible to have $\dim R <n$, for example when $w=x_1^{N+1}$ and $n\geq 2$, since in this case $\partial_{x_1} w= (N+1) x_1^N$ and 
$\partial_{x_i}w=0$ for $i>1$. Remark that if $\dim R <n$, then the complex $C_w$ is not isomorphic to its dual $C_w^{\vee}$ for a dimensional 
reason, therefore our complex $C_w$ is different from the complex considered by Bocklandt, Schedler and Wemyss 
(Lemma 6.4 in~\cite{bsw:higher}). According to the previous proposition, the assumption $\dim R=n$ is natural as far as the self-duality of $C_w$ is concerned. 
\Bcr \label{corprop}
Let $k$ be a field of characteristic zero. Let $V$ be an $n$-dimensional space with $n\geq 1$. Let $w$ be a non-zero homogeneous potential of $V$ of 
degree $N+1$ with $N\geq 2$. Let $A=A(w)$ be the potential algebra defined by $w$, so that the space of generators of $A$ is $V$. 
Assume that the space of relations $R$ of $A$ is $n$-dimensional. Then $f_0$, $f_1$, $f_2$ and $f_3$ form an isomorphism of complexes of graded $A$-bimodules 
$f: C_w \rightarrow C_w^{\vee}$, which is homogeneous of degree $-N-1$. Moreover, we have $\emph{H}^2 (A, A\stackrel{o}\otimes A)=0$. If the global 
dimension of $A$ is equal to 3, then $\emph{H}^0 (A, A\stackrel{o}\otimes A)=0$ and $\emph{H}^3 (A, A\stackrel{o}\otimes A)$ surjects onto $A$.
\Ecr
\Bdm
The first statement is immediate from Proposition \ref{commut}. Let us denote by $\bar{d}_3 : AEA \rightarrow ARA$ the first arrow in (\ref{start2}), 
while we keep $d_3: Akc(w)A \rightarrow ARA$ for the first arrow in (\ref{sdc}). So $d_3=\bar{d}_3 \circ i$ where $i$ denotes an obvious inclusion, and we 
have $d_3^{\ast}=i^{\ast}\circ \bar{d}_3^{\ast} $, hence $\ker \bar{d}_3^{\ast}$ is contained in $\ker d_3^{\ast}$. But the dual complex of (\ref{start2}) 
shows that $\mbox{im}\, d_2^{\ast}$ is contained in $\ker \bar{d}_3^{\ast}$. Carrying on the exactness of $C_w$ at $AV\!A$ by $f$, we finally get that 
$\mbox{im}\, d_2^{\ast}=\ker \bar{d}_3^{\ast}$, i.e. $\mbox{H}^2 (A, A\stackrel{o}\otimes A)=0$. The last statement is clear.
\Edm
\Bte \label{main2}
Let $k$ be a field of characteristic zero. Let $V$ be an $n$-dimensional space with $n\geq 1$. Let $w$ be a non-zero homogeneous potential of $V$ of 
degree $N+1$ with $N\geq 2$. Let $A=A(w)$ be the potential algebra defined by $w$, so that the space of generators of $A$ is $V$. 
If the space of relations $R$ of $A$ is $n$-dimensional, the following are equivalent.

(i) $A$ is 3-Calabi-Yau.

(ii) $A$ is AS-Gorenstein of global dimension 3.

(iii) $A$ is $N$-Koszul of global dimension 3 and $\dim R_{N+1}=1$.

(iv) The complex $C_w$ (see (\ref{sdc})) is exact in positive degrees.
\Ete
\Bdm
(i)$\Rightarrow$(ii) comes from Proposition 4.3 in~\cite{rbrt:pbw}. (ii)$\Rightarrow$(iii) follows from Proposition \ref{Gore}. 
(iii)$\Rightarrow$(iv) is obvious. If $C_w$ is exact in positive degrees, then it coincides with the Koszul resolution of $A$ which is self-dual by $f$ 
(Corollary \ref{corprop}), 
hence $\mbox{H}^i(A, A\stackrel{o}\otimes A) = 0$ whenever $i\neq 3$ and $\mbox{H}^3(A, A\stackrel{o}\otimes A)$ is isomorphic to $A$ as $A$-bimodule. 
Thus $A$ is 3-Calabi-Yau. This argument of self-duality was used by Bocklandt~\cite{bock:graded} (see also~\cite{rbrt:pbw}).
\Edm
\\

We are now ready to prove Theorem \ref{main} of the Introduction. 
\Bte \label{main3}
Let $k$ be a field of characteristic zero. Let $V$ be an $n$-dimensional space with $n\geq 1$. Let $w$ be a non-zero homogeneous potential of $V$ of 
degree $N+1$ with $N\geq 2$. Let $A=A(w)$ be the potential algebra defined by $w$. Assume that $A$ is $N$-Koszul. Then $A$ is 3-Calabi-Yau 
if and only the Hilbert series of the graded algebra $A$ is given by 
$$h_A(t)=(1-nt + nt^N - t^{N+1})^{-1}.$$
\Ete
\Bdm
 Suppose that 
 $$h_A(t)=(1-nt + nt^N - t^{N+1})^{-1}.$$
 Recall that the complex (\ref{start2}):
\begin{equation} 
AEA \stackrel{d_3}\longrightarrow A R A \stackrel{d_2}\longrightarrow
A V\! A \stackrel{d_1}\longrightarrow A A \rightarrow 0,
\end{equation}
is the beginning of a minimal projective resolution of $A$ in $A$-grMod-$A$. Since $A$ is $N$-Koszul, the assumption on $h_A(t)$ implies that 
$\dim R= n$, $\dim E=1$, and the global dimension of $A$ is equal to 3. In particular, one has $E=R_{N+1}=kc(w)$. Thus $A$ is 3-Calabi-Yau 
by Theorem \ref{main2}.

The converse is immediate from Proposition \ref{Gore} since any 3-Calabi-Yau algebra is AS-Gorenstein of global dimension 3. 
\Edm
\\

We will examine in Section 4 the recent examples which have motivated us to state this theorem. Now we want to show how this theorem allows us to 
recover some important examples of 3-Calabi-Yau homogeneous potentials. 
\Bex \Eex Let $S_n$ be the symmetric group of $\{ 1, \ldots , n\}$ and let $sgn$ be the sign of a permutation. Suppose that $n\geq 3$ is odd. Then 
the potential
$$c(w)=\mbox{Ant} (x_1, \ldots x_n) = \sum_{\sigma \in S_n} sgn(\sigma)\, x_{\sigma(1)},\ldots, x_{\sigma(n)}$$
is 3-Calabi-Yau. In fact, it is known~\cite{rb:nonquad} that 
$A(w)= T(V)/ I(\Lambda ^{n-1} V)$ is 
($n-1$)-Koszul of global dimension 3, $\dim R=n$ and $\dim R_{n} = 1$. Note that $w=\mbox{Ant}(x_1, \ldots, x_{n-1}) x_n$. The algebra 
$A(w)$ is called an antisymmetrizer algebra or an $(n-1)$-symmetric algebra in $n$ variables. It coincides with the polynomial algebra when $n=3$. If 
$n\geq 4$ is even, then $A(w)$ is not derived from a potential~\cite{rbrt:pbw}, but $A(w)$ can be derived from a super-potential defined by a super-cyclic sum $c$ 
given by the formula (1.2) in~\cite{bsw:higher}.
\Bex \Eex 
Assume that $n\geq 2$ and that $V$ is endowed with a non-degenerate symmetric bilinear form $(\,,\,)$. Setting $g_{ij}=(x_i,x_j)$, the inverse matrix of 
$(g_{ij})_{\leq i,j \leq n}$ is denoted by $(g^{ij})_{\leq i,j \leq n}$. Kriegk and Van den Bergh have shown in~\cite{kvdb:ncqg} that the space 
$Pot(V)_4^{O(V)}$ of the 4-degree potentials invariant by the orthogonal group $O(V)$ is 2-dimensional, generated by the following
$$w_1= \sum_{1\leq i,j,p,q\leq n} g^{ip} g^{jq}\, [x_i, x_j] [x_p,x_q],$$
$$w_2= ( \sum_{1\leq i,j\leq n} g^{ij} x_i x_j )^2,$$
where $[a,b]$ denotes the commutator of $a$ and $b$ in $T(V)$. For any $\lambda \in k$ with $\lambda \neq \frac{n-1}{n+1}$, the potential 
$w=w_1 + \lambda w_2$ is 3-Calabi-Yau. In fact, Kriegk and Van den Bergh have proved that, for any $\lambda \in k$, $A(w)$ is 3-Koszul and 
$\dim R =n$. For $\lambda \neq \frac{n-1}{n+1}$, Connes and Dubois-Violette have proved that $\dim R_4=1$ and $R_5=0$, which implies that the global 
dimension of $A(w)$ is equal to 3. The algebras $A(w)$ are called deformed Yang-Mills algebras (one omits ``deformed'' if $\lambda =0$) and were 
introduced by Connes and Dubois-Violette~\cite{cdv:ym, cdv:dym}.
\Bex \Eex
Fix $k=\mathbb{C}$ and three generators $x$, $y$, $z$. Set
$$S= \{(\alpha:\beta :\gamma)\in \mathbb{P}^2(\mathbb{C}) ; \alpha^3=\beta^3= 27 \gamma^3\}\cup \{(1:0:0), (0:1:0), (0:0:1)\}.$$
For any $(\alpha:\beta :\gamma)\in \mathbb{P}^2(\mathbb{C}) \setminus S$, the potential
$$w= \alpha xyz + \beta yxz + \gamma (x^3 + y^3 + z^3)$$
is 3-Calabi-Yau. In fact, the algebras $A(w)$ are exactly the generic quadratic AS-regular algebras of global dimension 3 and of type A (also called 
Sklyanin algebras in three generators), and one deduces from~\cite{as:regular} that $A(w)$ is Koszul with $h_A(t)=(1-t)^{-3}$.
\Bex \Eex
Fix $k=\mathbb{C}$ and two generators $x$, $y$. Set
$$S= \{(\alpha:\beta :\gamma)\in \mathbb{P}^2(\mathbb{C}) ; \alpha^2=4 \beta^2= 16 \gamma^2\}\cup \{(0:1:0), (0:0:1)\}.$$
For any $(\alpha:\beta :\gamma)\in \mathbb{P}^2(\mathbb{C}) \setminus S$, the potential
$$w= \alpha x^2y^2 + \beta (xy)^2 + \gamma (x^4 + y^4)$$
is 3-Calabi-Yau. In fact, the algebras $A(w)$ are exactly the generic cubic AS-regular algebras of global dimension 3 and of type A, and one deduces 
from~\cite{as:regular, rb:nonquad} that $A(w)$ is 3-Koszul with $h_A(t)=(1-2t+2t^3 -t^4)^{-1}$.
\Bex \label{ex5} \Eex
Consider $V$ of dimension 1, $V=k x$ and $w=x^{N+1}$. Then, $\dim R = \dim R_{N+1} =1$ and $A(w)$ is $N$-Koszul, but the global dimension of $A(w)$ is 
infinite. So $w$ is not 3-Calabi-Yau. Here is a non-trivial way to get other examples of $w$ which are not 3-Calabi-Yau. According to Theorem \ref{main2}, 
it suffices to assume that $\dim R =n$, $A(w)$ is $N$-Koszul of global dimension 3 and $\dim R_{N+1}>1$. Question: find such a potential.

\setcounter{equation}{0}

\section{Van den Bergh's duality}
From now on $k$ will be a field of characteristic zero. Let $V$ be a vector space of dimension $n\geq 2$, $w$ is a non-zero homogeneous potential 
of $V$ of degree $N+1$ with $N\geq 2$. Assume that the algebra $A=A(w)$ is 3-Calabi-Yau. Then, since $A$ is $N$-Koszul and AS-Gorenstein, 
it satisfies Van den Bergh's duality (Theorem 6.3 in~\cite{rbnm:kogo}), and the dualizing bimodule is $A$ itself. This means that for any $A$-bimodule 
$M$, there are linear isomorphisms between Hochschild homology and cohomology: $\mbox{H}_{\bullet}(A, M) \cong \mbox{H}^{3-\bullet}(A, M)$. 
We are going to construct an explicit isomorphism of complexes giving the above duality, from the self-duality 
$f: C_w \rightarrow C_w^{\vee}$ of the previous section. We do not use Van den Bergh's duality theorem~\cite{vdb:dual}.
 
Replace the assumption that $A$ is 3-Calabi-Yau by the weaker assumption $\dim R=n$ with $n\geq 1$, so that the self-duality $f$ 
still holds according to Corollary \ref{corprop}. From the chain complex isomorphism $f$, we define for any $A$-bimodule $M$ an isomorphism of complexes 
of $k$-vector spaces
\begin{equation} \label{iso}
M\otimes_{A^e} f : M\otimes_{A^e} C_w \longrightarrow M\otimes_{A^e} C_w^{\vee},
\end{equation}
and since $f$ has an inverse morphism $g$, then $M\otimes_{A^e} f$ has an inverse morphism which is $M\otimes_{A^e} g$. The flip 
$$\tau : M\otimes_{A^e} C_w^{\vee} \longrightarrow C_w^{\vee} \otimes_{A^e} M$$
is an isomorphism of complexes. Since $C_w$ is a projective left $A^e$-module of finite type, one has canonical isomorphisms of $k$-vector spaces 
(\cite{nb:alg}, Ch. 2, Prop. 2, p.75):
\begin{equation} \label{isocan}
C_w^{\vee} \otimes_{A^e} M = \mbox{Hom}_{A^e}(C_w, A^e) \otimes _{A^e} M \cong \mbox{Hom}_{A^e} (C_w, A^e \otimes_{A^e}M)
\cong \mbox{Hom}_{A^e} (C_w, M),
\end{equation}
so that the homology of the complex $M\otimes_{A^e} C_w^{\vee}$ is the Hochschild cohomology $\mbox{H}^{\bullet} (A,M)$. Thus, if $A$ is 3-Calabi-Yau, 
the isomorphism (\ref{iso}) in homology provides the expected linear isomorphisms $\mbox{H}_{\bullet}(A, M) \cong \mbox{H}^{3-\bullet}(A, M)$. Remark that 
if $A$ is not 3-Calabi-Yau, the complex $M\otimes_{A^e} C_w$ (resp. $M\otimes_{A^e} C_w^{\vee}$) only computes $\mbox{H}_0(A, M)$ and $\mbox{H}_1(A, M)$
(resp. $\mbox{H}^0 (A, M)$ and $\mbox{H}^1(A, M)$).

Let us express explicitly the isomorphism (\ref{iso}). We keep the weaker assumption $\dim R=n$ with $n\geq 1$. As previously, we omit the unadorned symbols 
$\otimes$ when they separate spaces. Moreover, $M\otimes kc(w)$ is denoted by $Mc(w)$. Firstly $M\otimes_{A^e} C_w$ is written down
\begin{equation}\label{comphom}
0 \longrightarrow Mc(w) \stackrel{\tilde d_3}{\longrightarrow}  MR\stackrel{\tilde d_2}{\longrightarrow} MV \stackrel{\tilde d_1}{\longrightarrow} M 
 \longrightarrow 0,
\end{equation}
where $M\otimes_{A^e} d$ is denoted by $\tilde d$. For any $m\in M$ and $1\leq i \leq n$, one has 
\begin{equation} \label{d1hom}
\tilde {d_1} (m\otimes x_i)= mx_i - x_i m = [m, x_i].
\end{equation}
From (\ref{d2}), we get 
$$\tilde{d_2} (m\otimes r_i)=\sum_{1\leq j\leq n} \sum_{1,2} \left(\frac{\partial r_i}{\partial x_j}\right)_2 m 
\left(\frac{\partial r_i}{\partial x_j}\right)_1 \otimes x_j .$$
Recalling that the entries of the Hessian matrix are given by
$$\frac{\partial^2 w}{\partial x_i \partial x_j}=\sum_{1,2} \left(\frac{\partial r_j}{\partial x_i}\right)_1 \otimes  
\left(\frac{\partial r_j}{\partial x_i}\right)_2=\sum_{1,2} \left(\frac{\partial r_i}{\partial x_j}\right)_2 \otimes  
\left(\frac{\partial r_i}{\partial x_j}\right)_1,$$
we see that 
\begin{equation} \label{d2hom}
\tilde{d_2} (m\otimes r_i)=\sum_{1\leq j\leq n} \left(\frac{\partial^2 w}{\partial x_i \partial x_j}\cdot m\right) \otimes x_j .
\end{equation}
From (\ref{eul2}), we get 
\begin{equation} \label{d3hom}
\tilde {d_3} (m\otimes c(w))= \sum_{1\leq i \leq n} [m, x_i]\otimes r_i.
\end{equation}

Next $M\otimes_{A^e} C_w^{\vee}$ is written down
\begin{equation}\label{comphomast}
0 \longrightarrow M \stackrel{\tilde {d_1^{\ast}}}{\longrightarrow}  MV^{\ast} \stackrel{\tilde {d_2^{\ast}}}{\longrightarrow} MR^{\ast} 
\stackrel{\tilde {d_3^{\ast}}}{\longrightarrow} Mc(w)^{\ast} \longrightarrow 0,
\end{equation}
where $M\otimes_{A^e} d^{\ast}$ is denoted by $\tilde {d^{\ast}}$. From (\ref{d1ast}), (\ref{d2ast}) and (\ref{d3ast}), we easily obtain
\begin{equation} \label{d1asthom}
\tilde {d_1^{\ast}} (m)= \sum_{1\leq i \leq n} [m, x_i]\otimes x_i^{\ast},
\end{equation}
\begin{equation} \label{d2asthom}
\tilde{d_2^{\ast}} (m\otimes x_i^{\ast})=\sum_{1\leq j\leq n} \left(\frac{\partial^2 w}{\partial x_i \partial x_j}\cdot m\right) \otimes r_j^{\ast},
\end{equation}
\begin{equation} \label{d3asthom}
\tilde {d_3^{\ast}} (m\otimes r_i^{\ast})= [m, x_i]\otimes c(w)^{\ast}.
\end{equation}

Finally, after applying the functor $M\otimes_{A^e}-$, the commutative diagram (\ref{self}) becomes 
\begin{eqnarray} \label{selfhom}
0 \rightarrow Mc(w) \stackrel{\tilde{d_{3}}}{\longrightarrow} & MR
\stackrel{\tilde{d_{2}}}{\longrightarrow} 
& MV \stackrel{\tilde{d_{1}}}{\longrightarrow} \ \ M \longrightarrow 0
\nonumber \\
\tilde{f_{3}} \downarrow \ \ \ \ \ \ \  & \tilde{f_{2}} \downarrow \ \ \ \ \  &  \tilde{f_{1}}
\downarrow \ \ \ \ \ \ \tilde{f_{0}} \downarrow \\
0\longrightarrow M \stackrel{\tilde{d_{1}^{\ast}}}{\longrightarrow} & \ \ MV^{\ast} \stackrel{\tilde{d_{2}^{\ast}}}{\longrightarrow}
& MR^{\ast} \stackrel{\tilde{d_{3}^{\ast}}}{\longrightarrow}  Mc(w)^{\ast} \rightarrow 0
\nonumber 
\end{eqnarray}
where $M\otimes_{A^e} f$ is denoted by $\tilde {f}$. One has immediately 
\begin{equation}
\tilde{f_0}(m)=m\otimes c(w)^{\ast}, \ \ \tilde{f_1}(m\otimes x_i)= m \otimes r_i^{\ast},\ \ 
\tilde{f_2}(m\otimes r_i)=m\otimes x_i^{\ast}, \ \ \tilde{f_3}(m\otimes c(w))= m.
\end{equation}

\setcounter{equation}{0}

\section{Skew polynomial algebras over n.-c. quadrics}
Let $u=\sum_{1\leq i,j \leq n} u _{ij} x_i x_j$ denote a quadratic polynomial in the 
non-commutative one-degree variables $x_1, \ldots , x_n$ with $n\geq 2$. We assume that $u$ is non-degenerate, meaning that the matrix $(u_{ij})_{1\leq i,j \leq n}$ is 
invertible. Let $\Gamma$ be the non-commutative quadric defined by $u$, i.e. $\Gamma$ is defined by
$$\Gamma = k\langle x_1, \ldots , x_n \rangle /I(u)$$ 
where $k\langle x_1, \ldots , x_n \rangle $ denotes the free associative algebra in $x_1, \ldots x_n$ and $I(u)$ the two-sided ideal generated by $u$. 
Then the graded algebra $\Gamma$ is Koszul, AS-Gorenstein of global dimension $2$, and $\Gamma$ is 2-Calabi-Yau if and only if $u$ is 
skew-symmetric~\cite{rb:gera}. Moreover, $\Gamma$ is left (right) noetherian if and only if $n=2$, and $\Gamma$ is always a domain~\cite{jz:nnoeth}. 
Its Hilbert series is given by $h_{\Gamma}(t)=(1-nt+t^2)^{-1}$.
Following~\cite{rbap:cyp}, let $z$ be an extra generator of degree 1, and let $w$ be a non-zero cubic homogeneous potential in the variables 
$x_1, \ldots x_n,z$. In the next proposition, $w$ is not necessarily equal to $uz$ as in~\cite{rbap:cyp}. Actually, we want to include in the same 
proposition the examples coming from~\cite{smith:octo, msa:steiner}. 
In~\cite{msa:steiner}, Su\'arez Alvarez treats the case of algebras coming from Steiner triple systems, generalizing the example studied in~\cite{smith:octo}. Let us explain briefly the form of the potentials considered in~\cite{msa:steiner}. Fix an oriented Steiner triple system $(E; S)$ of order $n+1$, where $E=\{1, \dots, n+1\}$ and $S$ is a set of directed triples in $E$ such that any pair in $E$ is included in a unique (non-directed) triple belonging to $S$. A directed triple is denoted by $(i,j,k)$, so that $(j,k,i)$ and $(k,i,j)$ are also directed (see~\cite{msa:steiner} for details on the choice of the orientation of the triples). Then Su\'arez Alvarez associates to $(E; S)$ the following potential 
\[w= \sum_{(i,j,k) \in S} (x_ix_j-x_jx_i)x_k , \]
where $x_{n+1}=z$. 

\smallskip

Let us go back to the generic case of a non-zero cubic homogeneous potential $w$ in the variables $x_1, \ldots x_n,z$. Let $A$ be the quadratic graded algebra derived from the potential $w$. 
\Bpo \label{apply}
We keep the above notation and assumptions. Assume that the graded algebra $A$ is isomorphic to a skew polynomial algebra $\Gamma [z;\sigma; \delta]$, 
where $\sigma$ is a 0-degree homogeneous automorphism of $\Gamma$ and $\delta$ is a 1-degree homogeneous $\sigma$-derivation of $\Gamma$. 
Then $A$ is Koszul and 3-Calabi-Yau, and it is a domain. Moreover, $A$ is left (right) noetherian if and only if $n=2$. 
\Epo
\Bdm
For a sketch on skew polynomial algebras, the reader is referred to~\cite{bg:aqg}, pp. 8-9. From $A\cong \Gamma[z;\sigma;\delta]$, we see that 
$h_A(t)= h_{\Gamma}(t)/1-t$, hence 
$$h_A(t)=(1-(n+1)t+(n+1)t^2-t^3)^{-1}.$$
Moreover, since $\Gamma$ is Koszul, then $A\cong \Gamma[z;\sigma;\delta]$ is Koszul (see Theorem 10.2 in~\cite{cs:K2}). 
Thus $A$ is 3-Calabi-Yau from Theorem \ref{main3}. 
Since $\Gamma$ is a domain, $A\cong \Gamma[z;\sigma;\delta]$ is a domain. If $\Gamma$ is left (right) noetherian, the same holds for 
$A\cong \Gamma[z;\sigma;\delta]$. If $A$ is left (right) noetherian, then the Gelfand-Kirillov dimension of $A$ is finite~\cite{sz:growth}, hence the poles 
of $h_A$ are complex numbers of module 1, implying that $n=2$. \Edm
\\ 

In the algebras considered in~\cite{rbap:cyp, smith:octo, msa:steiner}, the assumptions of Proposition \ref{apply} are satisfied (see the cited articles). 
Thus we recover that these algebras are Koszul and 3-Calabi-Yau. Note that $\delta =0$ in~\cite{rbap:cyp}, while $\sigma =Id_{\Gamma}$ and 
$\delta \neq 0$ in~\cite{smith:octo, msa:steiner}.
\\

For the rest of this section, we focus on the situation considered in~\cite{rbap:cyp}, that is, we assume that 
$$w=uz = \sum_{1\leq i,j \leq n} u_{ij}z.$$
Our aim is to show that the automorphism $\sigma$ in this case is related to Van den Bergh's duality of $\Gamma$. Denote by $V$ (resp. $V_{\Gamma}$) the 
space of generators of $A$ (resp. $\Gamma$), and by $R$ (resp. $R_{\Gamma}$) the corresponding space of relations. A basis of $V_{\Gamma}$ consists of 
$x_1, \dots, x_n$, and $V=V_{\Gamma} \oplus kz$. A basis of $R_{\Gamma}$ consists of $\partial_z w=u$, and it suffices to add 
$r_1=\partial_{x_1}w, \ldots , r_n=\partial_{x_n}w$ to get a basis of $R$. The graded algebra $\Gamma$ is isomorphic to the subalgebra of $A$ generated by 
$x_1, \ldots, x_n$ and it is also isomorphic to $A/I(z)$. Recall from~\cite{rbap:cyp} that the automorphism $\sigma$ of $A$ is defined for 
$1 \leq i \leq n$ by 
$$zx_i=\sigma (x_i)z.$$

For $1\leq i \leq n$, one has 
$$r_i= \sum_{1\leq j \leq n} (u_{ij} x_j z + u_{ji} z x_j).$$
Then it is easy to compute the entries of the Hessian matrix 
$$H=\left(\frac{\partial^2 w}{\partial x_i \partial x_j}\right)_{1\leq i,j \leq n+1}= 
\left(\frac{\partial r_j}{\partial x_i} \right)_{1\leq i,j \leq n+1}$$
for our potential $w$, where $x_{n+1}=z$ and $r_{n+1}=u$. For $1\leq i, j \leq n$, we find that
\begin{equation} \label{hess1}
\frac{\partial r_j}{\partial x_i} = u_{ji} 1\otimes z + u_{ij} z\otimes 1,
\end{equation}
\begin{equation} \label{hess2}
\frac{\partial r_j}{\partial z} = \sum_{1\leq i\leq n} (u_{ij} 1\otimes x_i + u_{ji} x_i\otimes 1),
\end{equation}
\begin{equation} \label{hess3}
\frac{\partial u}{\partial x_i} = \sum_{1\leq j\leq n} (u_{ij} 1\otimes x_j + u_{ji} x_j\otimes 1),
\end{equation}
\begin{equation} \label{hess4}
\frac{\partial u}{\partial z} = 0.
\end{equation}

The bimodule Koszul resolution $C_w$ of $A$ is given by (\ref{sdc}), that is
$$0 \rightarrow Akc(w)A \stackrel{d_3}\longrightarrow A R A \stackrel{d_2}\longrightarrow
A V\! A \stackrel{d_1}\longrightarrow AkA \rightarrow 0.$$
Let us apply the natural projections $A\rightarrow \Gamma$, $V\rightarrow V_{\Gamma}$ and $R\rightarrow R_{\Gamma}$ to this complex, in order to obtain 
the commutative diagram
\begin{eqnarray} \label{proj}
0 \rightarrow A(kc(w))A \stackrel{d_{3}}{\longrightarrow} & ARA
\stackrel{d_{2}}{\longrightarrow} 
& AV\! A \stackrel{d_{1}}{\longrightarrow} AkA \longrightarrow 0
\nonumber \\
 \downarrow \ \ \ \ \ \ \  \ \  &  \downarrow \ \ \ \ \ \ \ &  
\ \ \downarrow \ \ \ \ \ \ \ \ \ \ \ \   \downarrow  \\
0 \ \ \longrightarrow \ \ 0 \ \ \ \longrightarrow & \Gamma R_{\Gamma} \Gamma \stackrel{d_2^{\Gamma}}\longrightarrow
& \Gamma V_{\Gamma} \Gamma \stackrel{d_1^{\Gamma}}\longrightarrow \Gamma k \Gamma \rightarrow 0
\nonumber 
\end{eqnarray}
in which the second row is defined as follows. Since $c(w)$ vanishes 
modulo $z$, $d_3$ factors out to $0\rightarrow \Gamma R_{\Gamma} \Gamma$. Clearly, $d_1$ factors out to a $\Gamma$-$\Gamma$-linear map 
$d_1^{\Gamma}: \Gamma V_{\Gamma} \Gamma \rightarrow \Gamma \Gamma$ defined by $d_1^{\Gamma} (x_i)= x_i \otimes 1 -1 \otimes x_i$, $1\leq i \leq n$. 
From (\ref{d2}) and the entries of $H$ given above, one obtains for $1\leq j \leq n$, that 
$$d_2(r_j)=\sum_{1\leq i \leq n} u_{ij}(z\otimes x_i \otimes 1 + 1\otimes z \otimes x_i) + 
\sum_{1\leq i \leq n} u_{ji}(1\otimes x_i \otimes z +x_i \otimes z \otimes 1).$$
Thus $d_2$ factors out to a $\Gamma$-$\Gamma$-linear map $d_2^{\Gamma}: \Gamma R_{\Gamma} \Gamma \rightarrow \Gamma V_{\Gamma} \Gamma$. Using 
(\ref{hess3}) and (\ref{hess4}), $d_2^{\Gamma}$ is defined by 
\begin{equation} \label{d2gamma}
d_2^{\Gamma} (u)= \sum_{1\leq i,j \leq n} u_{ij} (1\otimes x_i \otimes x_j + x_i \otimes x_j \otimes 1).
\end{equation}
Consequently, the so-obtained quotient complex 
\begin{equation} \label{kgamma}
0 \rightarrow \Gamma R_{\Gamma} \Gamma \stackrel{d_2^{\Gamma}}\longrightarrow
\Gamma V_{\Gamma} \Gamma \stackrel{d_1^{\Gamma}}\longrightarrow \Gamma k \Gamma \rightarrow 0
\end{equation}
coincides with the bimodule Koszul resolution of $\Gamma$.

Let us proceed similarly with the dual complex $C_w^{\vee}$ given by (\ref{dsdc}):
$$0 \rightarrow Ak^{\ast}A \stackrel{d_1^{\ast}}\longrightarrow A V^{\ast} A \stackrel{d_2^{\ast}}\longrightarrow
A R^{\ast} A \stackrel{d_3^{\ast}}\longrightarrow Akc(w)^{\ast}A \rightarrow 0.$$
The natural inclusion $V_{\Gamma} \rightarrow V$ (resp. $R_{\Gamma} \rightarrow R$) provides the projection $V^{\ast} \rightarrow V_{\Gamma}^{\ast}$ (resp. 
$R^{\ast} \rightarrow R_{\Gamma}^{\ast}$). The image of the dual basis $(x_1^{\ast}, \ldots, x_n^{\ast}, z^{\ast})$ (resp. 
$(r_1^{\ast}, \ldots, r_n^{\ast}, u^{\ast})$) by this projection consists of 0 and the basis $(x_1^{\ast}, \ldots, x_n^{\ast})$ of $V_{\Gamma}^{\ast}$ 
(resp. the basis $u^{\ast}$ of $R_{\Gamma}^{\ast}$). Clearly, $d_3^{\ast}$ factors out to $\Gamma R_{\Gamma}^{\ast} \Gamma \rightarrow 0$ and 
$d_1^{\ast}$ factors out to $d_1^{\ast \Gamma}: \Gamma k^{\ast} \Gamma \rightarrow \Gamma V_{\Gamma}^{\ast} \Gamma $ defined by 
$$d_1^{\ast \Gamma}(1^{\ast})= \sum_{1\leq i \leq n} x_i \otimes x_i^{\ast} \otimes 1 - 1 \otimes x_i^{\ast} \otimes x_i.$$
From (\ref{d2ast}) and the entries of $H$, one deduces that
$$d_2^{\ast}(z^{\ast})= \sum_{1\leq i,j \leq n} u_{ij} (1\otimes r_j^{\ast} \otimes x_i + x_j \otimes r_i^{\ast} \otimes 1).$$
Since the RHS does not contain the element $u ^{\ast}$, $d_2^{\ast}$ factors out to 
$d_2^{\ast \Gamma}: \Gamma V_{\Gamma}^{\ast} \Gamma \rightarrow \Gamma R_{\Gamma}^{\ast} \Gamma $. Using (\ref{d2ast}) and (\ref{hess3}), 
$d_2^{\ast \Gamma}$ is defined for $1\leq i \leq n$ by 
\begin{equation} \label{d2astgamma}
d_2^{\ast \Gamma}(x_i^{\ast}) = \sum_{1\leq j \leq n} (u_{ij} x_j \otimes u^{\ast} \otimes 1 + u_{ji} 1 \otimes u^{\ast} \otimes x_j).
\end{equation}
Then, it is easy to show that the complex 
\begin{equation} \label{dkgamma}
0 \rightarrow \Gamma k^{\ast} \Gamma \stackrel{d_1^{\ast \Gamma}}\longrightarrow \Gamma V_{\Gamma}^{\ast} \Gamma \stackrel{d_2^{\ast \Gamma}}\longrightarrow
\Gamma R_{\Gamma}^{\ast} \Gamma \rightarrow 0
\end{equation}
is isomorphic to the dual complex of the complex of bimodules (\ref{kgamma}). In fact, following along the same lines of Section 2, it suffices to 
apply the isomorphism $\theta$ to this dual complex and to verify that we obtain the complex (\ref{dkgamma}). The verification is left to the reader.

Since $\Gamma$ is Koszul and AS-Gorenstein, $\Gamma$ satisfies Van den Bergh's duality (Proposition 2 in~\cite{vdb:dual}). More precisely, there is an 
automorphism $\nu$ of the graded algebra $\Gamma$ such that, for any $\Gamma$-bimodule $M$, we have linear isomorphisms 
$$\mbox{H}^{\bullet}(\Gamma, M) \cong \mbox{H}_{3-\bullet}(\Gamma,\, _{\nu}M).$$
As usual, the bimodule $_{\nu}M$ coincides with $M$ as right module but the left action of $a\in \Gamma$ upon $m\in \, _{\nu}M$ is given by $\nu(a)m$. In~\cite{vdb:dual}, $\nu$ is expressed 
in terms of the Nakayama automorphism of the dual Koszul algebra $\Gamma ^!$ (note that the same results hold for 
AS-Gorenstein $N$-Koszul algebras~\cite{rbnm:kogo}). In our situation, the following description of $\nu$ does not need the use of $\Gamma ^!$. 
\Bpo \label{autom}
We have $\nu= \sigma ^{-1}$. In particular, $\sigma =Id_{\Gamma}$ if and only if $u$ is skew-symmetric.
\Epo
\Bdm
Since $\Gamma$ satisfies Van den Bergh's duality, the homology of the complex (\ref{dkgamma}) at $\Gamma k^{\ast} \Gamma$ and at 
$\Gamma V_{\Gamma}^{\ast} \Gamma$ vanishes, and it is isomorphic to $_{\nu} \Gamma$ at $\Gamma R_{\Gamma}^{\ast} \Gamma$. Define the $\Gamma$-$\Gamma$-linear 
map $\mu_u : \Gamma R_{\Gamma}^{\ast} \Gamma \rightarrow \Gamma_{\sigma}$ by 
$$\mu_u (a\otimes u^{\ast} \otimes b)= a\, \sigma (b)$$
for any $a$ and $b$ in $\Gamma$. Choose $1\in \Gamma_{\sigma}$ of degree $-2$, so that $\mu_u$ is homogeneous of degree 0. Let us check that 
$\mu_u \circ d_2^{\ast \Gamma}=0$. Fix $i\in \{1, \ldots , n \}$ and set $X_i= \mu_u \circ d_2^{\ast \Gamma} (1\otimes x_i^{\ast} \otimes 1)$. From 
(\ref{d2astgamma}), we get that 
$$X_i= \sum_{1\leq j \leq n} (u_{ij} x_j + u_{ji} \sigma(x_j)).$$
Therefore $X_iz=\sum_{1\leq j \leq n} (u_{ij} x_j z + u_{ji} z x_j)=r_i$, hence $X_iz=0$ in $A$. But $z$ is not a zero-divisor in $A$, thus $X_i=0$ as 
desired. 

Next, examining the surjective homogeneous natural map
$$_{\nu} \Gamma  \cong \frac{\Gamma R_{\Gamma}^{\ast} \Gamma}{\mbox{im}d_2^{\ast \Gamma}} \rightarrow 
\frac{\Gamma R_{\Gamma}^{\ast} \Gamma}{\ker \mu_u} \cong \Gamma _{\sigma}$$
degree by degree, we see that it is an isomorphism and that $_{\nu} \Gamma  \cong \Gamma _{\sigma}$. Thus $\nu= \sigma ^{-1}$. In particular, 
$\sigma =Id_{\Gamma}$ if and only if $\Gamma$ is Calabi-Yau, i.e. if and only if $u$ is skew-symmetric, recovering 2) of Proposition 2.11 
in~\cite{rbap:cyp}.
\Edm
\\

\vspace{0.5 cm} 
\textsf{Roland Berger: Universit\'{e} de Lyon, Institut Camille Jordan (UMR 5208), Universit\'e de Saint-Etienne, 
Facult\'{e} des Sciences, 23, Rue Docteur Paul Michelon,   
42023 Saint-Etienne Cedex 2, France}

\emph{Roland.Berger@univ-st-etienne.fr}\\ 

\textsf{Andrea Solotar: IMAS y Dto de Matem\`{a}tica, Facultad de Ciencias Exactas y Naturales,
Universidad de Buenos Aires, 
Ciudad Universitaria, Pabell\`{o}n 1,
(1428) Buenos Aires, Argentina}

\emph{asolotar@dm.uba.ar}

\end{document}